\documentclass[12pt]{article}
\usepackage{latexsym,amsfonts,amssymb}
\usepackage{indentfirst}

\def\y{\begin{eqnarray*}}
\def\bd{\begin{description}}
\def\ey{\end{eqnarray*}}
\def\ebd{\end{description}}

\def\R{\mathbb{R}}

\begin{document}

\setlength{\baselineskip}{18pt}

\title{ Monotone iterative technique for  delayed evolution equation periodic problems in Banach spaces\footnote{Research supported by
NNSFs of China (11261053, 11361055).}}

\author{Qiang Li\footnote{Corresponding author.
E-mail: lznwnuliqiang@126.com (Q. Li)}\\[8pt]
\small{Department of Mathematics, Shanxi Normal University, }\\
\small{Linfen 041000, Peoples's Republic of China.}
}

\date{}

\maketitle

\begin{abstract}
\setlength{\baselineskip}{14pt}

In this paper, we deal with the existence of $\omega$-periodic mild solutions for the abstract evolution equation with delay in an ordered Banach space $E$
$$u'(t)+Au(t)=F(t,u(t),u(t-\tau)),\ \ \ \ t\in\R,$$
where $A:D(A)\subset E\rightarrow E$ is a closed linear operator and $-A$ generates a positive $C_{0}$-semigroup $T(t)(t\geq0)$,  $F:\R\times E\times E\rightarrow E$ is a continuous mapping which is $\omega$-periodic in $t$, and $\tau\geq0$ is a constant.
Under some weaker assumptions,
we construct   monotone iterative method for the delayed evolution equation periodic problems, and obtain the existence of maximal and minimal periodic mild solutions.
 The results obtained generalize the recent conclusions on this topic. Finally, we present
two applications to illustrate the feasibility of our abstract results.

\vspace{8pt}


\noindent {\bf  MR(2010) Subject Classification:\ } 34K30; 47H07; 47H08

\end{abstract}

\section{Introduction}

Let $E$ be an ordered Banach space, whose positive cone $K=\{u\in E | \ u\geq \theta\}$ is normal with normal constant $N$. In this paper, we use a monotone iterative technique in the presence of the lower and
upper solutions to discuss the existence of the extremal periodic mild solutions to the periodic
problem of first order semilinear evolution equations with delay in
 ordered Banach space $E$
$$u'(t)+Au(t)=F(t,u(t),u(t-\tau)),\qquad t\in\R, \eqno(1.1)$$
where $A:D(A)\subset E\rightarrow E$ is a closed linear operator and $-A$ generates a positive strongly
continuous semigroup ($C_{0}$-semigroup, in short)  $T(t)(t\geq0)$ in $E$, the nonlinear function $F:\R\times E\times E\rightarrow E$ is a continuous mapping  and for every $x,y\in E$,
 $F(t,x,y)$ is $\omega$-periodic in $t$ and $\tau$ is positive constant which denotes the time delay.

The theory of partial differential equations with delays
 is an important branch of differential equation theory, which has
 extensive physical, biological, economical,
  engineering background and realistic mathematical model,
   and hence has been emerging as an important area of
   investigation in the last few decades,
   and the numerous properties of their solutions
   have been studied, see \cite{Hale1993,Wu1996} and references therein.

The problem concerning periodic solutions of partial
 differential equations with delay is an important area
 of investigation since they can take into account
seasonal fluctuations occurring in the phenomena
appearing in the models, and have been studied by
some researchers in recent years. There has
been a significant development in periodic solution of evolution
equation with delay in Banach spaces,
 we refer to the references \cite{Burton1991,Xiang1992,Liu1998,Liu2000,Liu2003, Li2011,Kpoumie2016,LiLi2016,Wang2014,Liang2015,Liang2017}.

 In \cite{Burton1991}, Burton and Zhang obtained the existence
of periodic solutions for an abstract evolution equation
 with infinite delay. In \cite{Xiang1992}, under the assumption that
 the corresponding initial value problem has a priori estimate,
 Xiang and Ahmed showed an existence result of periodic solution to
the delay evolution equations in Banach spaces.
 In \cite{Liu1998,Liu2000,Liu2003}, Liu
derived periodic solutions from bounded solutions
or ultimate bounded solutions for finite or infinite
delay evolution equations in Banach spaces. In \cite{Li2011}, Li discussed the
periodic solutions of the evolution equation with delays and presented essential conditions on the nonlinearity to guarantee that the equation has periodic solutions.
In \cite{Kpoumie2016}, M. Kpoumi\`{e} et al. studied the existence of a periodic solution for some partial functional differential equations with infinite delay in Banach spaces. Recently, some authors also discussed the periodic solutions for some nonautonomous delay impulsive evolutionary equations (see \cite{Liang2015,Liang2017,Wang2014}). They established some existence results on periodic solutions to the equations under the ultimate boundedness of the solutions of the corresponding initial value problem.

In fact, in previous works, evolution equation periodic problems with delay have been studied
by many authors using different tools, such as Granas's fixed theorem, Banach contraction mapping principal, Schauder's fixed-point theorem, Horn's
fixed point theorem, Sadovskii's fixed point theorem and so on.
However, to the best of our knowledge,  few results yet exist for
the periodic problems with delay by using the method of the lower and upper solutions coupled
with the monotone iterative technique.

It is well known that the monotone iterative technique
of the lower and upper solutions
 is an effective and flexible mechanism.
 It yields monotone sequences of the lower
 and upper approximate solutions that converge to
 the minimal and maximal solutions between the lower
  and upper solutions. Early on, Du and Lakshmikantham \cite{DuLak82}, Sun and
Zhao \cite{SunZhao1992} investigated the existence of extremal solutions to the initial value problem
of ordinary differential equations without delay by using the method of the lower
and upper solutions coupled with the monotone iterative technique. Later, Li \cite{Li98} applied
 lower and upper solutions method to periodic
solution problems for semilinear evolution equations without delay in ordered Banach spaces, and obtained
the existence of maximal  and minimal periodic solutions using the characteristics of positive operators semigroups and the monotone iteration scheme.

Recently, in \cite{LiLi2015} we dealt with the second-order delayed ordinary differential equation periodic problem in ordered Banach spaces. With the nonlinear function satisfying quasi-monotonicity, we obtained the existence of the minimal and maximal periodic
solutions by monotone iterative technique of the lower and upper solutions. And in \cite{LiLi2016}, we also applied operator semigroup theory and monotone iterative technique of lower and upper solutions to obtain the existence and uniqueness of periodic mild solutions of
the abstract evolution equation under some quasi-monotone conditions.

Motivated by the papers mentioned above, the purpose of
 this paper is to construct the general principle
for lower and upper solutions coupled with the monotone iterative technique
 for the evolution equation periodic problems with delay, and obtain
the existence of maximal and minimal periodic mild solutions,
 which will make up the research in this area blank.

The paper is organized as follows. In Section 2,
some notions, definitions, and preliminary facts
 are introduced, which are used through this paper.
Under the different assumptions, the existence results
of the extremal periodic solutions of Equation (1.1)
 are given in Section 3. In Section 4, we  give two examples
  to illustrate our main  results in Section 3.

\section{Preliminaries}

In this section, we introduce some notions, definitions, and preliminary
facts which are used through this paper.

Throughout this paper, we assume that $E$ is an ordered Banach space, whose positive cone $K=\{u\in E|u\geq \theta\}$ is normal with normal constant $N$.

Let $A:D(A)\subset E\rightarrow E$ is a closed linear operator and  $-A$ generate a $C_{0}$-semigroup $T(t)(t\geq0)$ in $E$. For the theory of semigroups of linear operators we refer to \cite{Pazy83}.
We only recall here some notions and properties that are essential for us.
For a general $C_{0}$-semigroup $T(t)(t\geq0)$, there exist $M\geq1$ and $\nu\in \R$ such that (see \cite{Pazy83})
$$\|T(t)\|\leq Me^{\nu t},\quad t\geq0.\eqno(2.1)$$
Let
$$\nu_{0}=\inf\{\nu\in \R |\ \mathrm{There\ exists}\ M\geq1\ \mathrm{ such \ that}\  \|T(t)\|\leq Me^{\nu t},\ \forall t\geq0\},\eqno(2.2)$$
 then $\nu_{0}$ is called the growth exponent of the semigroup $T(t)(t\geq0)$. Furthermore, $\nu_{0}$ can be also obtained by the following formula
$$\nu_{0}=\limsup\limits_{t\rightarrow +\infty}\frac{\ln\|T(t)\|}{t}.$$
\vskip3mm
\noindent\textbf{Definition 2.1.}(\cite{Banasiak06}) \emph{ A $C_{0}$-semigroup $T(t)(t \geq 0)$ on $E$ is said to be positive, if the order inequality $T(t)x \geq \theta$ holds for each $x \geq\theta$, $x\in E$, and $t \geq 0$.}

It is easy to see that for any $C \geq 0$, $-(A + CI)$ also generates
 a $C_0$-semigroup $S(t) = e^{-Ct}T(t)(t\geq0)$ in $E$.
 And $S(t)(t\geq 0)$ is a positive $C_0$-semigroup if $T(t)(t \geq 0)$
  is a positive $C_0$-semigroup. For more details of the properties
   of the operator semigroups and positive $C_{0}$-semigroup,
we refer to the monographs \cite{Nagel86,Ruess94} and the paper \cite{Li1996}.

Let $J$ denote the infinite interval $[0,+\infty)$ and $h:J\rightarrow E$,
consider the initial value problem of the linear evolution equation
$$\left\{\begin{array}{ll}
u'(t)+Au(t)=h(t),\qquad t\in J,\\[8pt]
u(0)=x_{0}.
 \end{array} \right.\eqno(2.3)$$
 It is well known \cite[Chapter 4, Theorem 2,9]{Pazy83}, when $x_{0}\in D(A)$ and $h\in C^{1}(J,E)$,
 the initial value problem (2.3) has a unique classical solution
 $u\in C^{1}(J,E)\cap C(J,E_{1})$ expressed by
 $$u(t)=T(t)x_{0}+\int^{t}_{0}T(t-s)h(s)ds, \eqno(2.4)$$
where $E_{1}=D(A)$ is
Banach space  with the graph norm $\|\cdot\|_{1} = \|\cdot\|+\|A\cdot\|$. Generally, for $x_{0}\in E$ and $h\in C(J,E)$,
the function $u$ given by (2.4) belongs to $C(J,E)$ and it is
called a mild solution of the linear evolution equation (2.3).

Let $C_{\omega}(\R,E)$ denote the Banach space $\{u\in C(\R,E)|\ u(t)=u(t+\omega),t\in\R\}$ endowed the
maximum norm $\|u\|_{C}=\max_{t\in [0,\omega]}\|u(t)\|$. Evidently,
$C_{\omega}(\R,E)$ is  an order Banach space with the partial order $``\leq"$
induced by the positive cone $K_{C}=\{u\in C_{\omega}(\R,E)|\  \ u(t)\geq \theta,\ t\in \R\}$
and $K_{C}$ is also normal with the normal constant $N$. For $v, w \in C(\R,E)$ with $v \leq w$, we use $[v,w]$ to denote the order interval $\{u \in C(\R,E) | v \leq u \leq w\}$ in
$C(\R,E)$, and $[v(t), w(t)]$ to denote the order interval $\{u \in C(\R,E) | v(t) \leq u(t) \leq w(t), t \in \R\}$ in $E$.

Given $h\in C_{\omega}(\R,E)$, for the following linear evolution equation corresponding to Eq.(1.1)
$$u'(t)+Au(t)=h(t),\quad t\in \R,\eqno(2.5)$$
we have the following result.
\vskip3mm
\noindent\textbf{Lemma 2.2.}(\cite{Li2005x}) \emph{ If $-A$ generates an exponentially stable positive
$C_{0}$-semigroup $T(t)(t\geq0)$ in $E$, that is $\nu_{0}<0$ , then for $h\in C_{\omega}(\R,E)$,
the linear evolution equation (2.5) exists a unique positive $\omega$-periodic mild
solution $u$, which can be expressed by
  $$u(t)=(I-T(\omega))^{-1}\int_{t-\omega}^{t}T(t-s)h(s)ds:=(Ph)(t),\eqno(2.6)$$ and the solution operator $P:C_{\omega}(\R,E)\rightarrow C_{\omega}(\R,E)$ is a positive bounded linear operator.}
\vskip3mm
\noindent\textbf{Proof.} For any $\nu\in (0,|\nu_{0}|)$, there exists $M>0$ such that
$$\|T(t)\|\leq Me^{-\nu t}\leq M, \quad t\geq0.$$
In $E$, define the equivalent norm $|\cdot|$ by
$$|x|=\sup\limits_{t\geq0}\|e^{\nu t}T(t)x\|,$$
then $\|x\|\leq |x|\leq M\|x\|$. By $|T(t)|$ we denote the norm of $T(t)$
in $(E,|\cdot|$), then for $t\geq0$, it is easy to obtain that
$|T(t)|<e^{-\nu t}$.
Hence, $(I-T(\omega))$ has bounded inverse operator
$$(I-T(\omega))^{-1}=\sum_{n=0}^{\infty}T(n\omega),$$
and its norm satisfies
$$|(I-T(\omega))^{-1}|\leq \frac{1}{1-|T(\omega)|}\leq \frac{1}{1-e^{-\nu \omega}}.\eqno(2.7)$$
Set $$x_{0}=(I-T(\omega))^{-1}\int^{\omega}_{0}T(t-s)h(s)ds:=Bh,\eqno(2.8)$$
 then the mild solution $u(t)$  of  the linear initial value problem  (2.3)
given by (2.4) satisfies the periodic boundary condition $u(0)=u(\omega)=x_{0}$.
For $t\in \mathbb{R}^{+}$, by (2.4) and the properties of the semigroup $T(t)(t\geq0)$,
we have
\begin{eqnarray*}
u(t+\omega)&=&T(t+\omega)u(0)+\int_{0}^{t+\omega}T(t+\omega-s)h(s)ds\\
&=&T(t)\Big(T(\omega)u(0)+\int_{0}^{\omega}T(\omega-s)h(s)ds\Big)+\int^{t}_{0}T(t-s)h(s-\omega)ds\\
&=&T(t)u(0)+\int^{t}_{0}T(t-s)h(s)ds\ =u(t).
\end{eqnarray*}
Therefore, the $\omega$-periodic extension of $u$ on $\R$ is a unique $\omega$-periodic
mild solution of Eq.(2.5). By (2.4) and (2.8), the $\omega$-periodic mild solution
can be  expressed by
\begin{eqnarray*}\qquad\qquad\ \
u(t)&=&T(t)B(h)+\int^{t}_{0}T(t-s)h(s)ds\\
&=&(I-T(\omega))^{-1}\int^{t}_{t-\omega}T(t-s)h(s)ds:=(Ph)(t).\qquad\qquad\ \qquad (2.9)
\end{eqnarray*}
It is easy to see that  $P:C_{\omega}(\R,E)\rightarrow C_{\omega}(\R,E)$.
Finally, by the positivity
of semigroup $T(t)(t\geq0)$, we can obtain that $(I-T(\omega))^{-1}\geq\theta$, it follows that $Ph\geq\theta$ for any $h\in C_{\omega}(\R,E)$ and $h\geq\theta$. Therefore, $P:C_{\omega}(\R,E)\rightarrow C_{\omega}(\R,E)$ is a positive bounded linear operator.
This completes the proof of Lemma 2.2. $\Box$

Next, we recall some properties of measure of noncompactness that will
be used in the proof of our main results.
Let $\alpha(\cdot)$ denote the Kuratowski measure of noncompactness of the
bounded set. For the details of the definition and properties of the measure of noncompactness, see \cite{Banas80, Deimling85, Guo89}. For any $B\subset C_{\omega}(\R,E)$
and $t\in \R$, set $B(t) = \{u(t) | u \in B\} \subset E$. If $B$ is bounded in $C_{\omega}(\R, E)$, then $B(t)$ is bounded in $E$, and $\alpha(B(t))\leq \alpha(B)$.

 The following lemmas are needed in our arguments.
\vskip3mm
\noindent\textbf{Lemma 2.3.}(\cite{Banas80,Guo89,Guo06})\emph{ Let $E$ be a Banach space and let $B\subset C(J,E)$
be bounded and equicontinuous, where $J$ is a finite closed interval in  $\R$.
Then $\alpha(B(t))$ is continuous on $J$, and
$$ \alpha(B) = \max\limits_{t\in J}\alpha(B(t)) =\alpha(B(J)).$$}
\noindent\textbf{Lemma 2.4.}(\cite{Heinz83})\emph{ Let $E$ be a Banach space, $B=\{u_{n}\}\subset C(J,E) $
be a bounded and countable set. Then $\alpha(B(t))$ is Lebesgue
integrable on $J$, and
$$\alpha\Big(\Big\{\int_{J} u_n(s)ds\Big\}\Big)\leq
2\int_{J}\alpha(B(t))dt.$$}
\noindent\textbf{Lemma 2.5.} (\cite{Li01}) \emph{   Let $E$ be a Banach space and  $D\subset E$ be bounded. Assume that $Q:E\to E$ is linear bounded operator, then
$$\alpha(Q(D))\leq\|Q\|\alpha(D).
$$}
\section{Main results}

Now, we are in the  position to state and prove our main results. We will apply monotone iterative method of the lower and upper $\omega$-periodic solutions to obtain the existence of $\omega$-periodic mild solution for Eq.(1.1). To this end, we define the $\omega$-periodic lower and upper solutions of Eq.(1.1).

\vskip3mm
\noindent\textbf{Definition 3.1}  \emph{ If a function $v_0\in C^{1}_{\omega}(\R,E)\cap
C_{\omega}(\R,E_{1})$ satisfies
$$
v_0'(t)+Av_0(t)\leq F(t,v_0(t),v_{0}(t-\tau)),\;t\in \R,
\eqno(3.1)
$$
we call it an  $\omega$-periodic lower solution of Eq.(1.1). If  the inequality of (3.1) is inverse, we call it an  $\omega$-periodic upper solution of the Eq.(1.1).}

\vskip3mm
\noindent\textbf{Theorem 3.1}\quad \emph{Let $E$ be an ordered Banach space, whose positive cone $K$ is normal cone, let $A : D(A) \subset E \to E$ be
a closed linear operator and $-A$ generate a positive compact semigroup $T(t)(t\geq0)$, let $f:\R\times E\times E\rightarrow E$ be a continuous mapping which is $\omega$-periodic in $t$. Assume  Eq.(1.1) has lower and upper
$\omega$-periodic solutions $v_{0},w_{0}\in C_{\omega}^{1}(\R,E)\cap C_{\omega}(\R,E_{1})$ with $v_{0}\leq w_{0}$. If the following condition
\vskip1mm
\indent (H1)  there exists a constant $C\geq0$ such that for all $t\in\R$, $v_{0}(t)\leq x_{1}\leq x_{2}\leq w_{0}(t)$, $v_{0}(t-\tau)\leq y_{1}\leq y_{2}\leq w_{0}(t-\tau)$,
$$F(t,x_{2},y_{2})-F(t,x_{1},y_{1})\geq -C(x_{2}-x_{1})$$
\vskip1mm
\noindent holds, then the periodic problem (1.1)
has minimal and maximal $\omega$-periodic mild solution $\underline{u},\overline{u}$ between $v_{0}$ and $w_{0}$, which can
be obtained by monotone iterative sequences starting from $v_{0}$ and $w_{0}$. }
\vskip3mm
\noindent\textbf{Proof } Obviously, the periodic problem of  evolution equation with delay (1.1) is equal to the following  periodic problem
$$u'(t)+Au(t)+Cu(t)=F(t,u(t),u(t-\tau))+Cu(t),\ \ \ t\in\R, \eqno(3.2)$$
where the constant $C$ is decided by the condition (H1).

Let $C>|\nu_{0}|$(otherwise replace $C$ with $C+|\nu_{0}|$), then $-(A+CI)$ generates
 an exponentially stable, compact and positive $C_0$-semigroup $S(t) = e^{-Ct}T(t)(t\geq0)$ in $E$,
 whose growth exponent is $-C+\nu_{0}$. By Lemma 2.2, it follows that the following linear evolution equation periodic problem
$$u'(t)+Au(t)+Cu(t)=F(t,h(t),h(t-\tau))+Ch(t),\ \ \ t\in\R\eqno(3.3)$$
exists unique $\omega$-periodic mild solution
$$u(t)=(I-S(\omega))^{-1}\int_{t-\omega}^{t}S(t-s)(F(s,h(s),h(s-\tau))+Ch(s))ds:=Ph(t).\eqno(3.4)$$
From Definition 3.1, it is  clear that $[v_{0},w_{0}]\subset C_{\omega}(\R,E)$ and $v_{0}(t)\leq w_{0}(t)$ for any $t\in\R$.
 Define a mapping\ $\mathcal{F}:C_{\omega}(\R,E)\rightarrow C_{\omega}(\R,E)$ by
$$\mathcal{F}(u)(t)=F(t,u(t),u(t-\tau))+Cu(t),\ \ u\in C_{\omega}(\R,E), t\in\R,\eqno(3.5)$$
By the continuity of $F$, $\mathcal{F}:C_{\omega}(\R,E)\rightarrow C_{\omega}(\R,E)$ is continuous.
Define an operator $Q:[v_{0},w_0]\rightarrow C_{\omega}(\R,X)$ as follows:
$$Qu=(P\circ \mathcal{F})u,\ \ \eqno(3.6)$$
then we have
$$Qu(t)=(I-S(\omega))^{-1}\int_{t-\omega}^{t}S(t-s)(F(s,u(s),u(s-\tau))+Cu(s))ds,\ \ t\in \R,\eqno(3.7)$$
and $Q:[v_{0},w_{0}]\rightarrow C_{\omega}(\R,X)$ is continuous. Therefore, by the definition of $P$, we can assert $u\in [v_{0},w_{0}]$ is the $\omega$-periodic mild solution of Eq.(1.1) if and only if $u$ is the fixed point of the compound operator $Q$.

Now, we complete the proof by four steps.

\noindent\textbf{Step 1.} We show that the following properties of the operator $Q$ defined by (3.6).
\vskip0mm
\indent (i) $v_{0}\leq  Qv_{0}$ and\ $Qw_{0}\leq w_{0}$,
 \vskip0mm
\indent(ii) $Qu_{1}\leq Qu_{2}$ for any $u_{1},u_{2}\in[v_{0},w_{0}]$ with $u_{1}\leq u_{2}$.
\vskip0mm
\noindent Since\ $v_0\in C^{1}_{\omega}(\R,X)\cap
C_{\omega}(\R,X_{1})$ is an  $\omega$-periodic lower solution of Eq.(1.1),
thus
$$v_{0}'(t)+Av_{0}(t)+Cv_{0}(t)\leq F(t,v_{0}(t),
v_{0}(t-\tau))+Cv_{0}(t),\ \ \ t\in\R.\eqno(3.8)$$
Set\ $h(t)=v_{0}'(t)+Av_{0}(t)+Cv_{0}(t)$, by Lemma 2.2 and the positivity
of semigroup $S(t)(t\geq0)$, one can obtain that
$$v_{0}(t)=Ph(t)\leq P(F(t,v_{0}(t),v_{0}(t-\tau))+Cv_{0}(t))
\leq Qv_{0}(t),\ t\in\R,\eqno(3.9)$$
which implies that $v_{0}\leq Qv_{0}$. Similarly, it can be shown that $Qw_{0}\leq w_{0}$.

For any $u_{1},u_{2}\in[v_{0},w_{0}]$ with $u_{1}\leq u_{2}$ and $t\in\R$, we have
 $v_{0}(t)\leq u_{1}(t)\leq u_{2}(t)\leq w_{0}(t)$, $v_{0}(t-\tau)\leq u_{1}(t-\tau)\leq u_{2}(t-\tau)\leq w_{0}(t-\tau)$. By the condition (H1) and the positivity of the operator $P$,
\begin{eqnarray*}\qquad\quad \ \ Qu_{1}(t)&=&P(F(t,u_{1}(t),u_{1}(t-\tau))+C_{1}u(t))\\[8pt]
&\leq& P(F(t,u_{2}(t),u_{2}(t-\tau))+Cu_{2}(t))=Qu_{2}(t),\qquad\qquad\  (3.10)\end{eqnarray*}
 it follows that $Qu_{1}\leq Qu_{2}$.

 Therefore, $Q:
[v_0,w_0]\to [v_0,w_0]$ is a continuous increasing operator.

\noindent\textbf{Step 2.} We define two sequences $\{v_{i}\}$ and $\{w_{i}\}$ in $[v_{0},w_{0}]$ by the iterative scheme
$$v_{i}=Qv_{i-1},\quad w_{i}=Qw_{i-1},\ \qquad\qquad i=1,2,\cdots.\eqno(3.11)$$
Then from the monotonicity of the operator $Q$, it follows that
$$v_{0}\leq v_{1}\leq v_{2}\leq\cdots\leq v_{i}\leq\cdots\leq w_{i}
\leq\cdots\leq w_{2}\leq w_{1}\leq w_{0},\eqno(3.12)$$
and $\{v_{i}\},\{w_{i}\}\subset [v_{0},w_{0}]$ are equicontinuous in $\R$.

 In fact, for any $u\in [v_{0},w_{0}]$, by the periodicity of $u$, we consider it on $[0,\omega]$.
Set $ 0\leq t_{1}<t_{2}\leq\omega$, we get that
\begin{eqnarray*}
&&Qu(t_{2})-Qu(t_{1})\\[4pt]
&=&(I-T(\omega))^{-1}\int^{t_{2}}_{t_{2}-\omega}S(t_{2}-s)(F(s,u(s),u(s-\tau))+Cu(s))ds\\
&~&-(I-T(\omega))^{-1}\int^{t_{1}}_{t_{1}-\omega}S(t_{1}-s)(F(s,u(s),u(s-\tau))+Cu(s))ds\\
&=&(I-T(\omega))^{-1}\int^{t_{1}}_{t_{2}-\omega}(S(t_{2}-s)-S(t_{1}-s))(f(s,u(s),u(s-\tau))+Cu(s))ds\\
&~&-(I-T(\omega))^{-1}\int^{t_{2}-\omega}_{t_{1}-\omega}S(t_{1}-s)(f(s,u(s),u(s-\tau))+Cu(s))ds\\
&~&+(I-T(\omega))^{-1}\int^{t_{2}}_{t_{1}}S(t_{2}-s)(F(s,u(s),u(s-\tau))+Cu(s))ds\\
&:=&I_{1}+I_{2}+I_{3}.
\end{eqnarray*}
It is clear that
$$\|Qu(t_{2})-Qu(t_{1})\|\leq \|I_{1}\|+\|I_{2}\|+\|I_{3}\|.\eqno(3.13)$$
Thus, we only need to check $\|I_{i}\|$ tend to $0$  independently of $u\in[v_{0},w_{0}]$
when $t_{2}-t_{1}\rightarrow 0,i=1,2,3$.
For any $u\in[v_{0},w_{0}]$, from the condition (H1), it follows that \begin{eqnarray*}
F(t,v_{0}(t),v_{0}(t-\tau))+Cv_{0}(t)
&\leq& F(t,u(t),u(t-\tau))+Cu(t)\\[8pt]
&\leq& F(t,w_{0}(t),w_{0}(t-\tau))+Cw_{0}(t).
\end{eqnarray*}
By the normality of the cone $K$, there exists \ $M_{2}$ such that
$$\Big\|f(t,u(t),u(t-\tau))+Cu(t)\Big\|\leq M_{2}, \ \  t\in \R, u\in[v_{0},w_{0}].\eqno(3.14)$$
By the compactness of $S(t)(t\geq0)$, it follows that $S(t)$ is continuous in the uniform operator topology for $t > 0$. Hence, it is easy to check $\|I_{i}\|$ tend to $0$ independently of $u\in[v_{0},w_{0}]$
when $t_{2}-t_{1}\rightarrow 0(i=1,2,3)$,
which means that $Q([v_{0},w_{0}])$ is equicontinuous.

\noindent\textbf{Step 3.} $\{v_{i}(t)\}$ and  $\{w_{i}(t)\}$  are precompact on $E$ for any $t\in\R$.

Let $B_{1}=\{v_{i}\}$, $B_{2}=\{w_{i}\}$ and $B_{1}^{0}=B_{1}\cup\{v_{0}\}$, $B_{2}^{0}=B_{2}\cup\{w_{0}\}$.
Obviously  $t\in\R$, $B_{1}(t)=(QB_{1}^{0})(t)$  and  $B_{2}(t)=(QB_{2}^{0})(t)$ for $t\in\R$.

We define a set $(Q_{\varepsilon}B_{1})(t)$ by
$$(Q_{\varepsilon}B_{1}^{0})(t):=\{(Q_{\varepsilon}v_{i})(t)\
 |\ v_{i}\in B_{1}^{0},\ 0<\varepsilon<\omega,\ t\in\R\},\eqno(3.15)$$
where
\begin{eqnarray*}
&&Q_{\varepsilon}v_{i}(t)=(I-S(\omega))^{-1}
\int^{t-\varepsilon}_{t-\omega}S(t-s)\Big(f(s,v_{i-1}(t),v_{i-1}(s-\tau)
+Cv_{i}(s)\Big)ds\\[8pt]
&&~~=(I-S(\omega))^{-1} S(\varepsilon)\int^{t-\varepsilon}_{t-\omega}
S(t-s-\varepsilon)\Big(F(s,v_{i-1}(t),v_{i-1}(s-\tau)+Cv_{i}(s)\Big)ds.
\end{eqnarray*}
Then the set $(Q_{\varepsilon}B_{1}^{0})(t)$ is relatively compact in $E$ since the operator $S(\varepsilon)$ is compact in $E$ ($S(t)=e^{-Ct}T(t)(t\geq0) $ is compact semigroup). For any $v_{i}\in B_{1}^{0}$ and $t\in \R$, from the following inequality
\begin{eqnarray*}
&&\|Qv_{i}(t)-Q_{\varepsilon}v_{i}(t)\|\\[8pt]
&\leq&\Big\|(I-S(\omega))^{-1}\int_{t-\omega}^{t}S(t-s)
\Big(F(s,v_{i-1}(t),v_{i-1}(s-\tau))+Cv_{i-1}(s)\Big)ds\\[8pt]
&&-(I-S(\omega))^{-1}\int_{t-\omega}^{t-\varepsilon}S(t-s)
\Big(F(s,v_{i-1}(t),v_{i-1}(s-\tau))+Cv_{i-1}(s)\Big)ds\Big\|\\[8pt]
&\leq&\|(I-S(\omega))^{-1}\|\int_{t-\varepsilon}^{t}\Big\|S(t-s)
\Big(F(s,v_{i-1}(t),v_{i-1}(s-\tau))+Cv_{i-1}(s)\Big)\Big\|ds\\[8pt]
&\leq& \|(I-S(\omega))^{-1}\|M_{2}\int_{t-\varepsilon}^{t}\|S(t-s)\|ds,
\end{eqnarray*}
one can obtain that the set $(QB_{1}^{0})(t)$ is relatively compact, which implies that $\{v_{i}(t)\}=B_{1}(t)=(QB_{1}^{0})(t)$ is relatively compact in $E$ for $t\in \R$. Similarly, it can be shown that $\{w_{i}(t)\}$ is relatively compact in $E$ for $t\in \R$.

Therefore, $\{v_{i}\}$ and $\{w_{i}\}$ are relatively compact in $C_{\omega}(\R,E)$ by the Arzela-Ascoli Theorem, so there are convergent subsequences in $\{v_{i}\}$ and $\{w_{i}\}$, respectively. Combining this with the monotonicity and the normality of the cone $K_{C}$, we can easily prove that $\{v_{i}\}$ and $\{w_{i}\}$ themselves are convergent, i.e., there are $\underline{u},\overline{u}\in C_{\omega}(\R,E)$ such that $\lim\limits_{i\rightarrow\infty}v_{i}=\underline{u}$ and $\lim\limits_{i\rightarrow\infty}w_{i}=\overline{u}$.

Taking limit in (3.11), we have
$$\underline{u}=Q\underline{u},\qquad \overline{u}=Q\overline{u}.\eqno(3.16)$$
 Therefore $\underline{u},\overline{u}\in C_{\omega}(\R,X)$
are fixed points of $Q$, and they are the $\omega$-periodic mild solutions
of the periodic problem (1.1).

\textbf{Step 4}  We prove the minimal and maximal property of $\underline{u},\overline{u}$.

Assume that
$\widetilde{u}$ is a fixed point of $Q$ with $\widetilde{u}\in  [v_{0},w_{0}]$, then for every $t\in \R$, $v_{0}(t)\leq \widetilde{u}(t) \leq w _{0}(t)$,
$$v_{1}(t)= (Qv_{0})(t)\leq (Q\widetilde{u})(t)=\widetilde{u}(t)\leq (Qw_{0})(t)=w_{1}(t),\ t\in\R.\eqno(3.17)$$
Similarly, $v_{1}(t)\leq\widetilde{u}(t)\leq w_{1}(t)$, $t\in\R$.
In general
$$v_{i}\leq \widetilde{u}\leq w_{i},\ \quad  i=1,2,\cdots.\eqno(3.18)$$
Taking limit in (3.18) as $i\rightarrow\infty$, we get $\underline{u}\leq \widetilde{u}\leq\overline{ u}$.
Therefore $\underline{u},\overline{u}$ are minimal and maximal
 $\omega$-periodic mild solutions of Eq.(1.1), and $\underline{u},\overline{u}$
 can be obtained by the iterative sequences defined in (3.11) starting from
$v_{0}$ and $w_{0}$.
This completes the proof of Theorem 3.1.\hfill$\Box$
\vskip3mm
\noindent\textbf{Theorem 3.2} \emph{Let $E$ be an ordered Banach space, whose positive cone $K$ is normal cone, let $A : D(A) \subset E \to E$ be
a closed linear operator and $-A$ generate a positive equicontinuous $C_{0}$-semigroup $T(t)(t\geq0)$, let $F:\R\times E\times E\rightarrow E$ be a continuous mapping which is $\omega$-periodic in $t$. Assume the periodic problem (1.1) has lower and upper
$\omega$-periodic solutions $v_{0},w_{0}\in C_{\omega}^{1}(\R,E)\cap C_{\omega}(\R,E_{1})$ with $v_{0}\leq w_{0}$. If the condition (H1) and the following condition
\vskip0mm
\indent (H2) There exists a constant\ $c\in[0,1/4\omega C_{s}M_{S})$ such that for all\ $t\in\R$  and monotonic sequences\ $\{u_{n}\}\subset[v_{0},w_0]$,
$$\alpha(\{F(t,u_{n}(t),u_{n}(t-\tau))+Cu_{n}(t)\})\leq
 c(\alpha(\{u_{n}(t)\})+\alpha(\{u_{n}(t-\tau)\}))$$
hold, then the periodic problem (1.1) has minimal and maximal $\omega$-periodic mild solution $\underline{u},\overline{u}$ between $v_{0}$ and $w_{0}$, which can
be obtained by monotone iterative sequences starting from $v_{0}$ and $w_{0}$, where $C_{S}=\|(I-S(\omega))^{-1}\|$,
 $M_{S}=\sup\{\|S(t)\|\ |\  t\geq 0\}$.}
 \vskip3mm
\noindent\textbf{Proof } From the proof of Theorem 3.1, we know that $Q: [v_0,w_0] \to [v_0,w_0] $ is a continuous increasing operator and  $v_{0}\leq Qv_{0}$, $Qw_{0}\leq w_{0}$.
Hence, the iterative
sequences ${v_{i}}$ and ${w_{i}}$ defined by (3.11) satisfy (3.12).
By $T(t)(t\geq0)$ is an equicontinuous $C_{0}$-semigroup, it follows that $S(t)(t\geq0)$ is also an equicontinuous $C_{0}$-semigroup. From the proof of Theorem 3.1, we obtain that
$\{v_{i}\},\ \{w_{i}\}$ are equicontinuous in $\R$.

Next, we show that  $\{v_{i}\},\ \{w_{i}\}$  are
convergent in $C_{\omega}(\R,X)$.

Obviously, $\{v_{i}\}$ is a bounded countable set. By Lemma 2.4, Lemma 2.5 and the condition (H2), one can obtain that
\begin{eqnarray*}&&\alpha(\{v_{i}(t)\})=\alpha(\{Qv_{i-1}(t)\})\\[8pt]
&=&\alpha\Big(\Big\{(I-S(\omega))^{-1}\int_{t-\omega}^{t}S(t-s)
(F(s,v_{i-1}(s),v_{i-1}(s-\tau))+Cv_{i-1}(s))ds\Big\}\Big)\\[8pt]
&\leq&2\|(I-S(\omega))^{-1}\|\cdot\int_{t-\omega}^{t}\|S(t-s)\|
\cdot \alpha(\{F(s,v_{i-1}(s),v_{i-1}(s-\tau))+Cv_{i-1}(s)\})ds\\[8pt]
&\leq&2c\|(I-S(\omega))^{-1}\|\cdot\int_{t-\omega}^{t}\|S(t-s)\|
\cdot(\alpha(\{v_{i-1}(s)\})+\alpha(\{v_{i-1}(s-\tau)\})ds\\[8pt]
&\leq&2cC_{S}M_{S}\int_{t-\omega}^{t}\alpha(\{v_{i-1}(s)\})
+\alpha(\{v_{i-1}(s-\tau)\})ds,
\end{eqnarray*}
from the periodicity of $v_{i}$ and definition of measure of noncompactness, it follows that
 $\alpha(\{v_{i-1}(s)\})
=\alpha(\{v_{i-1}(s-\tau)\})$,
thus,
$$\alpha(\{v_{i}(t)\})\leq 4c\omega C_{S}M_{S}\cdot\max_{t\in[0,\omega]}\alpha(\{v_{i}(t)\})\leq 4c\omega C_{S}M_{S}\cdot\alpha_{C}(\{v_{i}\}).\eqno(3.19)$$
Since $\{v_{i}\}$ is equicontinuous, from Lemma 2.3, it follows that
$$0\leq \alpha_{C}(\{v_{i}\})\leq 4c\omega C_{S}M_{S}\cdot\alpha_{C}(\{v_{i}\}),$$
While $4c\omega C_{S}M_{S}<1$, hence\ $\alpha_{C}(\{v_{i}\})=0$.
Similarly, we can prove $\alpha_{C}(\{w_{i}\})=0$. Therefore, $\{v_{i}\},\ \{w_{i}\}$  are relatively compact in $C_{\omega}(\R,X)$,  so there are convergent subsequences in $\{v_{i}\}$ and $\{w_{i}\}$, respectively. Combining this with the monotonicity and the normality of the cone $K_{C}$, we can easily prove that $\{v_{i}\}$ and $\{w_{i}\}$ themselves are convergent, i.e., there are $\underline{u},\overline{u}\in C_{\omega}(\R,E)$ such that $\lim\limits_{i\rightarrow\infty}v_{i}=\underline{u}$ and $\lim\limits_{i\rightarrow\infty}w_{i}=\overline{u}$.

Therefore, from the proof of Theorem 3.1, $\underline{u},\overline{u}$ are minimal and maximal
 $\omega$-periodic mild solutions of the periodic problem with delay (1.1) in $[v_{0},w_{0}]$.
\hfill$\Box$
\vskip3mm

In the application of partial differential equations, we often choose Banach space $L^{p}(\Omega)(1\leq p<\infty)$
as working space, which is weakly sequentially complete space. Next, we discuss the existence of mild
solutions for the periodic problem with delay (1.1) in weakly sequentially complete Banach space.
\vskip2mm
\noindent\textbf{Theorem 3.3} \emph{ Let $E$ be an ordered and weakly sequentially complete Banach space, whose positive cone
$K$ is normal, let $A : D(A) \subset E \to E$ be
a closed linear operator and $-A$ generate a positive equicontinuous $C_{0}$-semigroup $T(t)(t\geq0)$ in $E$, let $F:\R\times E\times E\rightarrow E$ be a continuous mapping which is $\omega$-periodic in $t$.
Assume the periodic problem (1.1) has lower and upper
$\omega$-periodic solutions $v_{0},w_{0}\in C_{\omega}^{1}(\R,E)\cap C_{\omega}(\R,E_{1})$ with $v_{0}\leq w_{0}$. If the condition (H1) holds, then the periodic problem (3.1) has minimal and maximal $\omega$-periodic mild solution $\underline{u},\overline{u}$ between $v_{0}$ and $w_{0}$, which can
be obtained by monotone iterative sequences starting from $v_{0}$ and $w_{0}$.}
 \vskip3mm
\noindent\textbf{Proof }  From the proof of Theorem 3.1, it follows that the iterative
sequences ${v_{i}}$ and ${w_{i}}$ defined by (3.11) satisfy (3.12).
Hence, for any $t\in \R$, $\{v_i(t)\}$  and $\{w_i(t)\}$
are monotone and order-bounded sequences in $E$.
Noticing that $E$ is a weakly sequentially complete Banach
space, from Theorem 2.2 in  \cite{Du90}, one can get that $\{v_i(t)\}$ and$\{w_i(t)\}$  are precompact in $E$ for any $t\in \R$. Combining this
with the monotonicity (3.12), it follows that  $\{v_i(t)\}$ and$\{w_i(t)\}$ are uniformly convergent in $E$.
Denote
$$\underline{u}(t)=\lim\limits_{n\rightarrow\infty}v_n(t),\qquad
 \overline{u}(t)=\lim\limits_{n\rightarrow\infty}w_n(t),\;\;t\in \R. \eqno(3.20)$$
Obviously, $\{v_n(t)\}$, $\{w_n(t)\}\subset C_{\omega}(\R,X)$, and  $v_0(t)\leq
\underline{u}(t)\leq \overline{u}(t)\leq w_0(t)(t\in\R)$. By (3.7), we have
\begin{eqnarray*}
&&v_{i}(t)=Qv_{i-1}(t)\\[8pt]
&=&(I-S(\omega))^{-1}\int_{t-\omega}^{t}S(t-s)
(F(s,v_{i-1}(s),v_{i-1}(s-\tau))+Cv_{i-1}(s))ds
, \ \ (3.21)\end{eqnarray*}
and
\begin{eqnarray*}
&&w_{i}(t)=Qw_{i-1}(t)\\[8pt]
&=&(I-S(\omega))^{-1}\int_{t-\omega}^{t}S(t-s)
(F(s,w_{i-1}(s),w_{i-1}(s-\tau))+Cw_{i-1}(s))ds.  \ \ (3.22)\end{eqnarray*}
Taking limit in (3.21) and (3.22)  as $i\rightarrow\infty$, from the Lebesgue dominated
convergence theorem, one can obtain
$$
\underline{u}(t)=(I-S(\omega))^{-1}\int_{t-\omega}^{t}S(t-s)
(F(s,\underline{u}(s),\underline{u}(s-\tau))+C\underline{u}(s))ds
,\ \ t\in\R,\eqno(3.23)$$
and
$$
\overline{u}(t)=(I-S(\omega))^{-1}\int_{t-\omega}^{t}S(t-s)
(F(s,\overline{u}(s),\overline{u}(s-\tau))+C\overline{u}(s))ds
,\ \ t\in\R.\eqno(3.24)$$
which implies that $\underline{u},\overline{u}\in C_{\omega}(\R,X)$.
Similar with the proof of Theorem 3.1, we know that the $\underline{u},\overline{u}$ are minimal and maximal $\omega$-periodic mild solutions of the periodic problem with delay (1.1) in $[v_{0},w_{0}]$.
\hfill$\Box$

\vskip3mm
\noindent\textbf{Remark 1}\emph{ Analytic semigroup and differentiable semigroup are continuous by operator norm for every $t > 0$ (see \cite{Pazy83}). In the application of partial differential equations, such as parabolic equations and strongly damped wave equations, the corresponding solution semigroup is analytic semigroup. Therefore, Theorem 3.2 and Theorem 3.3 in this paper has broad applicability.}

In the above works, the key assumption (H1)
(the monotone on the third variable of the nonlinear function)
 is employed. However, we hope that the nonlinear function
  is quasi-monotonicity. In this case, the results have more extensive application background.

In fact, we find that if the periodic problem (1.1) has lower and upper
$\omega$-periodic solutions $v_{0},w_{0}\in C_{\omega}^{1}(\R,E)\cap C_{\omega}(\R,E_{1})$ with $v_{0}\leq w_{0}$ and
\vskip1mm
\noindent \emph{(H3) there is a sufficiently small constant $C_{1}>0$, such that
$$u_{2}(t)-u_{1}(t)\geq C_{1}( u_{2}(t-\tau)-u_{1}(t-\tau)), \ \ \ t\in\R,$$
for any $u_{1},u_{2}\in [v_{0},w_{0}]$ with\ $u_{2}\geq u_{1}$,}
\vskip1mm
\noindent then the condition (H1) can be replaced by the following condition
\vskip1mm
 \noindent  \emph{(H4) there are nonnegative constants\ $C_{2}, C_{3}$, such that
$$F(t,x_{2},y_{2})-F(t,x_{1},y_{1})
\geq -C_{2}(x_{2}-x_{1})-C_{3}(y_{2}-y_{1}),$$
for all $t\in \R$, $x_{1},x_{2},y_{1},y_{2}\in E$
with $v_{0}(t)\leq x_{1}\leq x_{2}\leq w_{0}(t)$,
\ $v_{0}(t-\tau)\leq y_{1}\leq y_{2}\leq w_{0}(t-\tau)$.}
\vskip1mm
 \indent  In fact, for every  $t\in\R$ and $u_{1},u_{2}\in [v_{0},w_{0}]$
 with $u_{1}\leq u_{2}$, one can obtain that  $v_{0}(t)\leq u_{1}(t)\leq u_{2}(t)\leq w_{0}(t)$,
 $v_{0}(t-\tau)\leq u_{1}(t-\tau)\leq u_{2}(t-\tau)\leq w_{0}(t-\tau)$.
 By  the conditions (H3) and (H4), it follows that
\begin{eqnarray*}
 &&F(t,u_{2}(t),u_{2}(t-\tau))-F(t,u_{1}(t),u_{1}(t-\tau))\\[8pt]
&\geq& -C_{2}(u_{2}(t)-u_{1}(t))-C_{3}(u_{2}(t-\tau)-u_{1}(t-\tau))\\[8pt]
&\geq& -C_{2}(u_{2}(t)-u_{1}(t))-\frac{C_{3}}{C_{1}}(u_{2}(t)-u_{1}(t))\\[8pt]
&=&-(C_{2}+\frac{C_{3}}{C_{1}})(u_{2}(t)-u_{1}(t))\\[8pt]
&:=&-C(u_{2}(t)-u_{1}(t)).
\end{eqnarray*}
\vskip1mm
\noindent Hence, we can obtain the following results form Theorem 3.1 and Theorem 3.2, respectively.
\vskip3mm
\noindent\textbf{Theorem 3.4}\quad \emph{Let $E$ be an ordered Banach space, whose positive cone $K$ is normal cone, let $A : D(A) \subset E \to E$ be
a closed linear operator and $-A$ generate a positive compact semigroup $T(t)(t\geq0)$, let $f:\R\times E\times E\rightarrow E$ be a continuous mapping which is $\omega$-periodic in $t$. Assume  Eq.(1.1) has lower and upper
$\omega$-periodic solutions $v_{0},w_{0}\in C_{\omega}^{1}(\R,E)\cap C_{\omega}(\R,E_{1})$ with $v_{0}\leq w_{0}$. If the  conditions (H3) and (H4) hold, then the periodic problem (1.1)
has minimal and maximal $\omega$-periodic mild solution $\underline{u},\overline{u}$ between $v_{0}$ and $w_{0}$, which can
be obtained by monotone iterative sequences starting from $v_{0}$ and $w_{0}$. }
\vskip3mm

\noindent\textbf{Theorem 3.5}\quad \emph{Let $E$ be an ordered Banach space, whose positive cone $K$ is normal cone, let $A : D(A) \subset E \to E$ be
a closed linear operator and $-A$ generate a positive equicontinuous $C_{0}$-semigroup $T(t)(t\geq0)$, let $f:\R\times E\times E\rightarrow E$ be a continuous mapping which is $\omega$-periodic in $t$. Assume the periodic problem (1.1) has lower and upper
$\omega$-periodic solutions $v_{0},w_{0}\in C_{\omega}^{1}(\R,E)\cap C_{\omega}(\R,E_{1})$ with $v_{0}\leq w_{0}$. If the  conditions (H2-H4) hold, then the periodic problem (1.1)
has minimal and maximal $\omega$-periodic mild solution $\underline{u},\overline{u}$ between $v_{0}$ and $w_{0}$, which can
be obtained by monotone iterative sequences starting from $v_{0}$ and $w_{0}$. }

\vskip3mm
\noindent\textbf{Remark 2}\emph{ Obviously, the condition (H3) is easy to satisfy, and the condition (H4) weakens the condition (H1). Hence, Theorem 3.3 and Theorem 3.4 partially improve Theorem 3.1 and Theorem 3.2. }

Next, we discuss the uniqueness of the $\omega$-periodic mild solution for the periodic problem (1.1) under $T(t)(t\geq0)$ is an equicontinuous $C_0$-semigroup.

\vskip3mm

\noindent\textbf{Theorem 3.6}\quad \emph{ Let $E$ be an ordered Banach space, whose positive cone $K$ is normal cone with normal constant $N$, let $A : D(A) \subset E \to E$ be
a closed linear operator and $-A$ generate a positive equicontinuous $C_{0}$-semigroup $T(t)(t\geq0)$, let $f:\R\times E\times E\rightarrow E$ be a continuous mapping which is $\omega$-periodic in $t$. Assume Eq.(1.1) has lower and upper
$\omega$-periodic solutions $v_{0},w_{0}\in C_{\omega}^{1}(\R,E)\cap C_{\omega}(\R,E_{1})$ with $v_{0}\leq w_{0}$. If the conditions (H3),(H4) and
   \vskip0mm
 \indent (H5) there exist constants  $L_{1},L_{2}>0$, such that for every $t\in \R$
and $x_{1},x_{2},y_{1},y_{2}\in X$, satisfying $v_{0}(t-\tau)\leq y_{1}\leq y_{2}\leq w_{0}(t-\tau)$,
$v_{0}(t)\leq x_{1}\leq x_{2}\leq w_{0}(t)$,
 $$F(t,x_{2},y_{2})-F(t,x_{1},y_{1})
\leq L_{1}(x_{2}-x_{1})+L_{2}(y_{2}-y_{1}),$$
\vskip2mm
 \indent (H6)  $N\Big(C_{2}+\frac{C_{3}}{C_{1}}+L_{1}+C_{1}L_{2}\Big)C_{S}M_{S}\omega<1$,
\vskip2mm
 \noindent hold, then the periodic problem (1.1)
has a unique $\omega$-periodic mild solution $u^{*}\in[v_{0},w_{0}]$, where  $C_{S}=\|(I-S(\omega))^{-1}\|$, $M_{S}=\sup\{\|S(t)\|\ |\ t\geq 0\}$.}
 \vskip3mm
\noindent\textbf{Proof }  From Theorem 3.1 and Theorem 3.5, one can obtain that
the iterative
sequences ${v_{i}}$ and ${w_{i}}$ defined by (3.11) satisfy (3.12).
 For any $t\in\R$,
   by the conditions (H3), (H5), (3.7), (3.12), it is to see
   \begin{eqnarray*}
   \theta&\leq&w_{i}(t)-v_{i}(t)=Qw_{i-1}(t)-Qv_{i-1}(t)\\[8pt]
   &=&(I-S(\omega))^{-1}\int_{t-\omega}^{t}S(t-s)
   (F(s,w_{i-1}(s),w_{i-1}(s-\tau))+C w_{i-1}(s))ds\\[8pt]
   &&-(I-S(\omega))^{-1}\int_{t-\omega}^{t}S(t-s)
   (F(s,v_{i-1}(s),v_{i-1}(s-\tau))+C v_{i-1}(s))ds\\[8pt]
   &\leq&(I-S(\omega))^{-1}\int_{t-\omega}^{t}S(t-s)
   ((L_{1}+C)(w_{i-1}(s)-v_{i-1}(s))\\[8pt]
   &&+L_{2}(w_{i-1}(s-\tau)-v_{i-1}(s-\tau)))ds\\[8pt]
   &\leq&(L_{1}+C+L_{2}C_{1})(I-S(\omega))^{-1}\int_{t-\omega}^{t}
   S(t-s)(w_{i-1}(s)-v_{i-1}(s))ds,\end{eqnarray*}
where $C=C_{2}+\frac{C_{3}}{C_{1}}$.
 By the normality of the cone $K$, it follows that
 $$\|w_{i}(t)-v_{i}(t)\|\leq
 N(L_{1}+C+L_{2}C_{1})C_{S}M_{S}\omega\|w_{i-1}-v_{i-1}\|_{C},\ \ t\in\R,$$
 namely
 $$\|w_{i}-v_{i}\|_{C}\leq
 N(L_{1}+C+L_{2}C_{1})C_{S}M_{S}\omega\|w_{i-1}-v_{i-1}\|_{C},\eqno(3.25)$$
 by the condition (H6), we can obtain that
 $$\|w_{i}-v_{i}\|_{C}\leq
 \Big(N(L_{1}+C+L_{2}C_{1})C_{S}M_{S}\omega\Big)^{i}\|w_{0}-v_{0}\|_{C}
 \rightarrow 0,\ \ i\rightarrow\infty.$$
 Thus, there is a unique $\omega$-periodic mild solution $u^\ast\in
C_{\omega}(\R,X)$, such that $\lim\limits_{i\rightarrow\infty}w_i
=\lim\limits_{i\rightarrow\infty}v_i=u^\ast$.
 Hence, taking limit in (3.11) as $i\rightarrow\infty$, we get $u^\ast=Qu^\ast$, which implies that $u^\ast$ is  unique $\omega$-periodic mild solution $u^{*}\in[v_{0},w_{0}]$ of the periodic problem (1.1).
\hfill$\Box$
\vskip3mm

\section{Application}

In this section, we present two examples, which do not aim at generality but indicate how our abstract
results can be applied to concrete problems.

\noindent\textbf{Example 4.1}\emph{ Periodic solutions of delay parabolic equations in $\R^{n}(n\geq1)$.}

Let $\overline{\Omega}\in\R^{n}$ be a bounded domain with a sufficiently smooth boundary $\partial\Omega$.
Let
$$A(x,D)u=-\sum^{N}_{i,j=1}a_{ij}(x)D_{i}D_{j}u+\sum_{j=1}^{N}a_{j}(x)D_{j}u+a_{0}(x)u,\eqno(4.1)$$
 be a uniformly elliptic differential operator in $\overline{\Omega}$, whose coefficients $a_{ij}(x), a_{j}(x)$ $(i,j=1,\cdots,n)$ and $a_{0}(x)$ are H\"{o}der-continuous on $\overline{\Omega}$, and $a_{0}(x)\geq0$.
 We let $B=B(x,D)$ be a
boundary operator on $\partial \Omega$ of the form:
$$Bu:=b_{0}(x)u+\delta\frac{\partial u}{\partial\beta},\eqno(4.2)$$
where either $\delta=0$ and $b_{0}(x)\equiv1$ (Dirichlet boundary operator), or $\delta=1$ and $b_{0}(x)\geq0$ (regular oblique derivative boundary operator; at this point, we further assume that $a_{0}(x)\not \equiv0$ or $b_{0}(x)\not \equiv0$), $\beta$ is an outward pointing, nowhere tangent vector field on $\partial\Omega$.
Let $\lambda_{1}$ be the first eigenvalue of elliptic operator  $A(x,D)$  under the boundary condition
 $Bu=0$. It is well known (\cite[Theorem 1.16]{Amann76},) that $\lambda_{1}>0$.

Under the above assumptions, we discuss the existence and uniqueness of
$\omega$-periodic solutions of the semilinear parabolic equation boundary value problem
$$\left\{\begin{array}{ll}
\frac{\partial }{\partial t}u(x,t)+A(x,D) u(x,t)=f(x,t,u(x,t),u(x,t-\tau)),\ x\in \Omega, \ t\in\R,\\[10pt]
Bu=0,\quad x\in \partial\Omega,
 \end{array} \right.\eqno (4.3)$$
where\ $f:\overline{\Omega}\times\R\times\R^{2}\rightarrow\R$ a local H\"{o}lder-continuous function which is $\omega$-periodic in $t$, $\tau>0$ denotes the time delay.

\vskip1mm
\noindent\textbf{Theorem 4.1}\emph{  Let $f:\overline{\Omega}\times\R\times\R^{2}\rightarrow\R$ be a local H\"{o}lder-continuous function which is $\omega$-periodic in $t$. If the following conditions
\vskip0mm
 \indent(H7) $f(x,t,0,0)\geq0$ for any $(x,t)\in\overline{\Omega}\times\R$, and there is a function $0\leq w=w(x,t)\in C^{2,1}(\Omega\times\R)$ which is $\omega$-periodic in $t$, such that
 $$\left\{\begin{array}{ll}
\frac{\partial }{\partial t}w(x,t)+A(x,D) w(x,t)\geq f(x,t,w(x,t),w(x,t-\tau)),(x,t)\in \Omega\times\R,\\[10pt]
Bw=0,\quad x\in \partial\Omega,
 \end{array} \right.$$
\vskip0mm
 \indent(H8) there exists a constant $c>0$, such that for any $x\in \Omega\ ,t\in \R$ and
 $0\leq y_{1}\leq y_{2}\leq w(x,t)$, $0\leq z_{1}\leq z_{2}\leq w(x,t-\tau)$,
$$f(x,t,y_{2},z_{2})-f(x,t,y_{1},z_{1})\geq-C(y_{2}-y_{1}),$$
\vskip0mm
\noindent hold, then the semilinear delayed parabolic equation boundary value problem (4.3) has minimal and maximal $\omega$-periodic  solution $\underline{u},\overline{u}\in  C^{2,1}(\overline{\Omega}\times\R)$ between $0$ and $w$, which can be obtained by monotone iterative sequences starting from $0$ and $w$.}
\vskip1mm
\noindent\textbf{Proof} Let\ $E=L^{p}(\Omega)(p>1)$, $K=\{u\in E|\ u(x)\geq0 \ a.e.\  x\in\Omega\}$,  then\ $E$ is an ordered Banach
space, whose positive cone $K $ is a normal regeneration cone.
Define an operator\ $A:D(A)\subset E\rightarrow E$ by:
$$D(A)= \{u\in W^{2,p}(\Omega)|\ B(x,D)u=0,\ x\in \partial\Omega\},\quad Au=A(x,D)u.\eqno(4.4)$$
If $a_{0}(x)\geq0$, then $-A$ generates an exponentially stable analytic semigroup $T_{p}(t)(t\geq0)$
in $E$ (see \cite{Amann78}). By the maximum principle of elliptic operators, we know that $(\lambda I+A)$ has a positive bounded inverse operator $(\lambda I+A)^{-1}$ for $\lambda>0$, hence $T_{p}(t)(t\geq0)$ is a positive
semigroup (see \cite{Li1996}). From the operator $A(x,D)$  has compact resolvent in $L^{p}(\Omega)$, we obtain $T_{p}(t)(t\geq0 )$ is also a compact semigroup (see \cite{Pazy83}).

Denote\ $u(t)=u(\cdot,t)$,  and $F(t,u(t),u(t-\tau))=f(\cdot,t,u(\cdot,t),u(\cdot,t-\tau))$, then parabolic boundary value problem (4.3) can be reformulated as the abstract evolution (1.1) in $E$.
By the condition (H7), it follows that $v_{0}\equiv0$ and $w_{0}=w(x,t)$ are time $\omega$-periodic lower solution and  time $\omega$-periodic upper solution of  the problem (4.3), and $v_{0}\leq w_{0}$. By the condition (H8), it follows that the condition (H1) holds. Hence, form Theorem 3.1, one can see the  delayed parabolic boundary value problem (4.3)
has minimal and maximal $\omega$-periodic mild solution $\underline{u},\overline{u}$, which can be obtained by monotone iterative sequences starting from $0$ and $w$, respectively.

By the analyticity
of the semigroup $T_{p}(t)(t\geq0)$ and the regularization method used in \cite{Amann78} , we can see that
$\underline{u},\overline{u}\in C^{2,1}(\overline{\Omega}\times\R)$  are   time $\omega$-periodic solutions of the problem(4.3). This completes
the proof of the theorem. \hfill$\Box$

 Furthermore, if the following condition
\vskip1mm
\indent \emph{(H9) there is a constant $C_{1}>0$, such that
$$u_{2}(x,t)-u_{1}(x,t)\geq C_{1}( u_{2}(x,t-\tau)-u_{1}(x,t-\tau)), $$
for $(x,t)\in\overline{\Omega}\times \R, u_{1},u_{2}\in [0,w(x,t)], u_{2}\geq u_{1}$,
}
then the condition (H8) can be replaced by
\vskip1mm
 \indent \emph{(H10) there exist nonnegative constants $C_{2}, C_{3}$ such that
 $$f(x,t,y_{2},z_{2})-f(x,t,y_{1},z_{1})
\geq -C_{2}(y_{2}-y_{1})-C_{3}(z_{2}-z_{1}).$$
for  $(x,t)\in \overline{\Omega}\times\R$
and  $x_{1},x_{2},y_{1},y_{2}\in X$
with $0\leq x_{1}\leq x_{2}\leq w_{0}(x,t)$,
\ $0\leq y_{1}\leq y_{2}\leq w_{0}(x,t-\tau)$.}

Thus, according to Theorem 3.4, we have the following result

\vskip1mm
\noindent\textbf{Theorem 4.2}\emph{  Let $f:\overline{\Omega}\times\R\times\R^{2}\rightarrow\R$ be a local H\"{o}lder-continuous function which is $\omega$-periodic in $t$. If the conditions (H7),(H9) and (H10) hold, then the semilinear delayed parabolic equation boundary value problem (4.3) has minimal and maximal $\omega$-periodic  solution $\underline{u},\overline{u}\in  C^{2,1}(\overline{\Omega}\times\R)$ between $0$ and $w$, which can be obtained by monotone iterative sequences starting from $0$ and $w$.}

\vskip2mm
\noindent\textbf{Example 4.2}\emph{ Doubly periodic problems of first order partial differential equation with delay.}

Let $f:\R^{4}\to \R$ is a continuous function, which is $2\pi$-periodic in $t$ and $x$. We are concerned with the existence of solutions for the semilinear first order partial differential equation with delay in $\R^{2}$:
$$\frac{\partial}{\partial t}u(x,t)+\frac{\partial}{\partial x}u(x,t)=f(x,t, u(x,t), u(x,t-\tau)),\ \ \  (x,t)\in\R^{2},\eqno(4.5)$$
with doubly periodic boundary conditions
$$u(x+2\pi,t)=u(x,t+2\pi)=u(x,t),\ \ \ (x,t)\in\R^{2},\eqno(4.6)$$
where $\tau>0$ denotes the time delay.

\noindent\textbf{Theorem 4.3}\emph{ Let $f(x,t,u,v)\in C^{1}(\R^{4})$, and $f$ is $2\pi$-periodic in $t$ and $x$. Assume $f(x,t,0,0)\geq 0$, there is a function $w(x,t)\in C^{1}(\R^{2})$ and $w$ is $2\pi$-periodic in $t$ and $x$ satisfying $w(x,t)\geq0$,  such that
 $$\frac{\partial}{\partial t}w(x,t)+\frac{\partial}{\partial x}w(x,t)\geq f(x,t, w(x,t), w(x,t-\tau)), \ \ \ \ (x,t)\in\R^{2}.$$
 If the following conditions
 \vskip1mm
 \indent(H11) for any $x,t\in \R$
and $u_{1},u_{2}\in C(\R^{2})$
, $0\leq u_{1}(x,t)\leq u_{2}(x,t)\leq w(x,t)$, \begin{eqnarray*}&f(x,t,u_{2}(x,t),u_{2}(x,t-\tau))-f(x,t,u_{1}(x,t),u_{1}(x,t-\tau))&\\[8pt]
&\geq -(u_{2}(x,t)-u_{1}(x,t)),& \end{eqnarray*}
\vskip-2mm
 \indent(H12) there is a constant\ $c\in [0,\frac{e^{2\pi}-1}{8\pi e^{2\pi}})$,
 such that for any  $x,t\in \R$  and monotone sequence\ $\{u_{n}(x,t)\}\in[0,w(x,t)]$,
  \begin{eqnarray*}&\alpha(\{f(x,t,u_{n}(x,t),u_{n}(x,t-\tau))+u_{n}(x,t)\})&\\[8pt]
  &\geq  c(\alpha(\{u_{n}(x,t)\})+\alpha(\{u_{n}(x,t-\tau)\})),&\end{eqnarray*}
 \vskip-2mm
 \noindent hold, then the doubly periodic problems of first order partial differential equation (4.5)-(4.6) has minimal and maximal  classical solutions $\underline{u},\overline{u}\in  C^{1}(\R^{2})$ between $0$ and $w$.}
\vskip1mm
\noindent\textbf{Proof} Let $C_{2\pi}(\R)$ denote the Banach space $\{u\in C(\R)| u(x+2\pi)=u(x), x\in \R\}$ endowed the maximum norm $\|u\|_{C}=\max_{x\in[0,2\pi]}\|u(x)\|$.
Denote $E=C_{2\pi}(\R)$, let
$$D(A)=C_{2\pi}^{1}(\R),\ \ A=\frac{\partial u}{\partial x} .\eqno(4.7)$$

From \cite[Lemma 2.1]{Li04}, if $\lambda\neq 0$, we know that $(\lambda I+A)$ has a bounded inverse operator $(\lambda I+A)^{-1}$ in $E$ and
$$(\lambda I+A)^{-1}h(x)=\int^{x}_{x-2\pi}r(s-y)h(y)dy, \ \ \ h\in E,\eqno(4.8)$$
where
$$r(x)=\frac{e^{-\lambda x}}{1-e^{-2\pi \lambda}},\ \ \ \ x\in[0,2\pi].$$
By\ (4.8), it follows that $(\lambda I+A)^{-1}$ is positive operator for $\lambda>0$, and its norm  $\|(\lambda I+A)^{-1}\|\leq \frac{1}{\lambda}$. Form  Hille-Yosida Theorem and exponential formula of semogroup(see \cite{Pazy83}), we can obtain that $-A$ generates a contractive and positive $C_{0}$-semigroup $T(t)(t\geq0)$, whose growth exponent $\nu_{0}\leq0$. Thus, $-(A+I)$  generates
 a contractive and positive $C_0$-semigroup $S(t) = e^{-t}T(t)(t\geq0)$ in $E$, and the growth exponent $\nu_{1}=-1+\nu_{0}\leq-1$, which implies that $S(t)(t\geq0)$ is an exponentially stable, positive $C_0$-semigroup and $\|S(2\pi)\|\leq e^{-2\pi}$, $\|(I-S(2\pi))^{-1}\|\leq\frac{e^{2\pi}}{e^{2\pi}-1}$.

Set $u(t)(x)=u(x,t),u(t-\tau)(x)=u(x,t-\tau)$, and
$$F(t,u(t),u(t-\tau))(x)=f(t, u(x,t), u(x,t-\tau)),\eqno(4.9)$$
then the doubly periodic problems (4.5)-(4.6) can be reformulated as following
$$u'(t)+Au(t)=F(t,u(t),u(t-\tau)),\ \ \ t\in\R,\eqno(4.10)$$
where $F:\R\times\ E\times E\to E$ is \ $C_{1}$-mapping which is $2\pi$-periodic in $t$.

It is easy to see that $v_{0}(t)\equiv0$ and $w_{0}(\cdot,t)=w(x,t)$ are  $2\pi$-periodic lower solution and   $2\pi$-periodic upper solution of Eq.(4.10). From the condition (H11),(H12) and Theorem 3.2,  one can obtain that Eq.(4.10) has minimal and maximal $2\pi$-periodic mild solution $\underline{u},\overline{u}$ between $0$ and $w_{0}$, which can
be obtained by monotone iterative sequences starting from $0$ and $w_{0}$. Since $F$ is a $C^{1}$-mapping, from regularity of solutions for the semilinear evolution equations (see\cite{Pazy83}), we know that
$$\underline{u},\overline{u}\in C_{2\pi}^{1}(\R,X)\cap  C_{2\pi}(\R,D(A)),\eqno(4.11)$$
namely $\underline{u},\overline{u}$ are minimal and maximal $2\pi$-periodic classical solutions, respectively.
 Therefore, by the definition of $A$, it follows that $\underline{u},\overline{u}$ are $2\pi$-doubly periodic  classical solutions of the doubly periodic problems  (4.5)-(4.6).
\hfill$\Box$

\bibliographystyle{abbrv}

\end{document}